\documentclass[12pt]{amsart}
\usepackage{amsfonts, amsmath, amsthm,amscd,amssymb,verbatim,epsf}
\usepackage{psfig}
\usepackage{amsmath}
\usepackage{amscd}
\newtheorem{pro}{Proposition}[section]
\newtheorem{thm}[pro]{Theorem}
\newtheorem{lem}[pro]{Lemma}

\newtheorem{cnj}[pro]{Conjecture}
\newtheorem{cor}[pro]{Corollary}

\theoremstyle{definition}

\theoremstyle{remark}

 \def\d{{\delta}}

 \def\l{{\lambda}}

   \def\s{{\sigma}}
 
 \def\a{{\alpha}}
 \def\b{{\beta}}
 \def\p{{\partial}}
 
 \def\ra{{\rightarrow}}

 \def\D{{\Delta}}

 \def\2{{\mathbb Z_2}}
 
 \def\t{{\tau}}
 
 \def\th{{\theta}}

\def\d{\delta}
\def\D{\Delta}
\def\a{\alpha}
\def\b{\beta}
\def\t{\tau}

\def\l{\lambda}

\def\s{\sigma}

\def\wh{\widehat}
\def\wt{\widetilde}

\begin{document}

\title{Virtually Haken surgeries on once-punctured torus bundles}
\author{Joseph D. Masters}
\begin{abstract}
 We describe a class $\mathcal{C}$ of punctured torus bundles such
 that, for each $M \in \mathcal{C}$,
 all but finitely many Dehn fillings on $M$ are virtually Haken.
 We show that $\mathcal{C}$ contains infinitely many commensurability
 classes, and we give evidence that $\mathcal{C}$ includes
 representatives of ``most'' commensurability classes
 of punctured torus bundles.

 In particular, we define an integer-valued complexity function
 on monodromies $f$ (essentially the length of the LR-factorization
 of $f_*$ in $PSL_2(\mathbb{Z})$),
 and use a computer to show that if the monodromy of $M$
 has complexity at most 5, then
 $M$ is finitely covered by an element of $\mathcal{C}$.
 If the monodromy has complexity at most 12,
 then, with at most 36 exceptions,
 $M$ is finitely covered by an element of $\mathcal{C}$.

 We also give a method for computing ``algebraic boundary slopes''
 in certain finite covers of punctured torus bundles.
\end{abstract}

\maketitle

\section{Introduction}

 A compact 3-manifold $M$ is \textit{Haken} if it is irreducible,
 and contains an orientable, essential surface. $M$ is
 \textit{virtually Haken} if it is finitely covered by a Haken manifold.
 One of the central problems in 3-dimensional topology
 is Waldhausen's Conjecture, which states that every closed, irreducible
 3-manifold with infinite fundamental group
 is virtually Haken.

 A \textit{knot manifold} is an irreducible, orientable, compact 3-manifold
 whose boundary is a single torus.  A knot manifold $M$
 is \textit{small} if every closed incompressible surface
 in $M$ is boundary parallel.
 A knot manifold is \textit{hyperbolic} if its interior admits a
 complete hyperbolic metric of finite volume.
 A \textit{slope} on a torus $T$ is a non-trivial
 isotopy class of simple closed
 curves on $T$. If $M$ is a 3-manifold,
 $T$ is a torus component of $\partial M$,
 and $\alpha$ is a slope on $T$, then $M(\alpha)$ denotes the manifold
 obtained from $M$ by Dehn filling along a simple closed curve
 representing the slope $\alpha$.
 We say that a knot manifold has \textit{Property VH}
 if $M(\alpha)$ is virtually Haken for all but finitely many
 slopes $\alpha$ on $\partial M$. 
 An important special case of Waldhausen's conjecture is the following:

\begin{cnj} \label{cnj}
Let $M$ be a hyperbolic knot manifold.  Then $M$ has Property VH.
\end{cnj}

 Examples of small knot manifolds satisfying Property VH have
 been given in \cite{BZ},\cite{M} and \cite{DT}.

 Let $F$ be homeomorphic to a torus with an open disk removed.
  We choose a basepoint, $p \in \partial F$, for $\pi_1 F$.
 Let $x$ and $y$ be generators of $\pi_1 F$, as pictured in Figure 1.
 Let $D_x$ and $D_y$ represent Dehn twists along simple closed curves
 in $int \, F$
 which are isotopic to loops representing $x$ and $y$ respectively.
 Let $\mathcal{M}_1^1$ denote the mapping class group of $F$;
 that is the group of orientation-preserving automorphisms of $F$
 which restrict to the identity on $\partial F$, modulo
 isotopies which fix every point of $\partial F$.
 It is well-known that $\mathcal{M}_1^1$ is generated by $D_x$ and
 $D_y$.  There is a natural map
 $\phi: \mathcal{M}_1^1 \rightarrow  SL_2(\mathbb{Z})$;
 we sometimes use the notation $\phi(f) = f_*$.

 Let $H_3 = \phi^{-1} < D_{x*}, D_{y*}^3>$,
 and let $H_4 = \phi^{-1} < D_{x*}, D_{y*}^4>$
 It is a fact (see Section 7) that $H_3$ and $H_4$
 are both finite-index subgroups of $\mathcal{M}_1^1$.

\begin{thm} \label{main}
Let $M$ be an orientable, atoroidal 3-manifold which fibers over $S^1$,
 whose fiber, $F$, is a compact, orientable surface of genus 1,
 with a single boundary component. Let $f:F \rightarrow F$
 be the monodromy, and
 suppose that the mapping class of $f$ lies in the subgroup $H_i$, where
 $i=3$ or $4$. Then there  are slopes
 $\beta_i^1, \beta_i^2$ and an integer $N$ such that
 $M(\alpha)$ is virtually Haken whenever
 $I(\alpha, \beta_i^1) > N$ and $I(\alpha, \beta_i^2) > 1$.
\end{thm}

Since $H_3$ and $H_4$ have finite index, it follows
 that every hyperbolic punctured torus bundle has infinitely many
 virtually Haken surgeries, a result which was first proved by Baker \cite{B}
 (see also \cite{M}). The slopes $\beta_i^j$ are computable (see Sections 6-8).

 For a given monodromy $f$, let $\b_i^j = \b_i^j(f)$ be slopes
 as in Theorem \ref{main}.
Let $\mathcal{C}$ be the set of mapping classes $f$ of $F$
 such that $f \in H_3 \cap H_4$ and
 $\{ \beta_3^1, \beta_3^2 \} \cap \{ \beta_4^1, \beta_4^2 \} = \emptyset$.

\begin{cor}
If $f \in \mathcal{C}$, then $M$ has Property VH.
\end{cor} 

 We shall show that there are infinitely many commensurability
 classes of once-punctured torus bundles whose
 monodromies lie in the class $\mathcal{C}$.  Thus we obtain:

\begin{thm} \label{inf}
There are infinitely many pair-wise non-commensurable
 once-punctured torus bundles which have Property VH.
\end{thm}

\noindent
\textit{Remarks:}\\
1. In the case where $f \in H_4$, Theorem \ref{main} is a corollary
 of Theorem 1.3 in \cite{M}.  The proof we give here is completely
 different.
\\
\\
2. Many of the bundles in $\mathcal{C}$ have no exceptional surgeries, and
 so the techniques of \cite{BZ} and \cite{DT} cannot be applied.\\

 It appears that (in some sense) ``most'' monodromies $f$
 have a power which lies in $\mathcal{C}$,
 and thus most punctured torus bundles are commensurable to ones
 with Property VH.
 We have verified this on a computer, for monodromies of low complexity. 

 To make these statements precise, we introduce a complexity
 function on monodromies.
 Recall that every element $g \in SL_2(\mathbb{Z})$
 can be written uniquely as  a positive word in $D_{x*}^{-1}$
 and $D_{y*}$, times $\pm Id$.
 We define the ``complexity'' of $g$ to
 be the length of this word.
 Similarly, if $f \in \mathcal{M}_1^1$, then we define the
 complexity of $f$ to be the complexity of $f_*$ in $SL_2(\mathbb{Z})$,
 and if $M$ is a punctured-torus bundle with monodromy
 $f$, the complexity of $M$ is defined
 to be the minimum of the complexity of $f$ among
 all monodromies of bundles $N$ which are bundle equivalent to $M$. 
 Using a computer, we can show:

\begin{thm} \label{computer}
a.  Every mapping class of complexity at most 5 has a power which
 lies in $\mathcal{C}$.
 Thus every once-punctured torus bundle of complexity at most 5
 is commensurable to one with Property VH.

b.There are 745 once-punctured torus bundles of complexity at most 12;
 with at most  36 exceptions, these
 are all commensurable to bundles with Property VH.\\
 \end{thm}

\noindent
\textit{IDEA OF PROOF}\\
 Our proof of Theorem \ref{main} is inspired by
 the arguments of Cooper and Walsh (\cite{CW}),
 who show that every every fibered knot in a $\mathbb{Z}/2$ homology
 sphere admits infinitely many virtually Haken surgeries.
 The idea is to replace $M$ with a finite cover $\widetilde{M}$
 in which the fiber has multiple boundary components.
 Then one hopes to find a non-separating surface in $\widetilde{M}$
 which is not a fiber, so that the techniques of \cite{CL}
 may be applied to find an
 essential surface in a cyclic cover of $\widetilde{M}$.
 Finally one must show that certain slopes on $\partial M$
 lift to the ultimate cover.

 Our argument diverges from that of \cite{CW}
 in that we choose $\widetilde{M}$ to have three or four boundary components
 instead of two.  Since we have two covers to work with,
 we obtain two surfaces, which is the key
 to proving Property VH.  Furthermore, the surfaces which
 we construct are disjoint from one of the boundary components
 of $\wt M$, and thus cannot be fibers.
 Thus we avoid a number of issues in \cite{CW} involving
 semi-bundle structures.

 However, we encounter several new issues.
 First, it is possible that our surfaces may become fibers after Dehn
 filling, and to rule this out requires the computation of
 an Alexander polynomial.
 Secondly, and more importantly, we require
 that the boundary components of the
 non-separating surface in $\widetilde{M}$
 must all project to the same slope in $M$. 
 To arrange this, we must develop techniques for
 constructing surfaces, and computing slopes, explicitly. 
 Some of these techniques (those in Section 3)
 may be applied to any bundle, but some (those in Section 4)
 exploit special features of the genus 1 mapping
 class group.\\
\\
\textit{A QUESTION}\\
 Although some of the methods in this paper apply only to punctured
 torus bundles, it is conceivable Conjecture \ref{cnj} can
 be attacked along broadly similar lines.  A key step would be to answer
 the following (presumably difficult) question:\\
\\
\textbf{Question:} Let $M$ be a knot manifold.  Is there a finite
 cover $\wt M$ of $M$ which contains a non-separating surface,
 which is disjoint from some component of $\partial \wt M$, and
 whose boundary curves all project to the same (embedded)
 slope on $\partial M$?\\
\\
By the results of this paper, the answer to the question is yes
 if $M$ is a punctured torus bundle.\\
\\
 \textit{PLAN OF PAPER}\\
Section 2 fixes a choice of basis for $H_1 (\partial M)$,
 when $M$ is any 3-manifold fibering over $S^1$.
 In Section 3, we show how to compute ``algebraic'' boundary slopes
 of non-separating surfaces in bundles.
 In Section 4, we prove the main theorem.
 In Sections 5-7, we give methods for explicit
 computations of slopes.
 Section 8 is devoted to the example
 of the figure eight knot exterior.  In Section 9, we prove Theorem \ref{inf}.
 Finally, in Section 10, we discuss our
 computer-generated data.\\
\\
\textit{ACKNOWLEDGMENT}\\
\\
Thanks are due to Genevieve Walsh for pointing out an error in a
 previous version.

\section{Framing convention}

 Let $f:F \rightarrow F$ be an automorphism of a compact,
 orientable surface, and let $M = F \times [0,1]/ (x,0) = (fx, 1)$
 be the 3-manifold fibering over $S^1$ with monodromy $f$.
 Let  $\lambda_1, ..., \lambda_k$ be the boundary components of $F$,
 and suppose that $f$ acts trivially
 on $\partial F$, so $M$ has torus boundary components
 $T_1, ..., T_k$.

 We wish to fix a framing for $\partial M$ (see Figure 1).
 For $i = 1, ..., k$,
 we fix a point $p_i \in \lambda_i$, and let the meridian,
 $\mu_i \subset T_i$, be the suspension of the point $p_i$.
 The orientation of $\mu_i$ is chosen so that
 the map from $[0,1]$ (with standard orientation) to  $\mu_i$
 given by  $t \rightarrow (p_i, 1-t)$ is orientation-preserving.
 We let the longitude, $\lambda_i$, of
 $T_i$ be given by $\lambda_i \times \{ 1 \}$.
 We orient $\lambda_i$ so that $I(\mu_i, \lambda_i) = 1$,
 where $I(.,.)$ is the standard intersection pairing on $H_1(\partial M)$.
 (see Figure 1).

 Given a surface $S$ properly embedded in $M$, we may specify
 the homology classes of the boundary curves of $S$ by a vector
 $(\alpha_1, ..., \alpha_k)$, where each $\alpha_i$
 is an ordered pair of integers.
 We shall refer to this as the vector of ``algebraic boundary slopes''
 of $S$.
 For example if we say that $S$ has algebraic boundary slopes
 $((1,2), (0,0), (0,3))$, we mean that
 $[S \cap T_1] = [\mu_1 + 2 \lambda_1]$ in $H_1(T_1)$,
 that $[S \cap T_2] = [0] \in H_1(T_2)$,
 and that $[S \cap T_3] =  3[\lambda_3]$ in $H_1(T_3)$.

\begin{figure}[!htbp]
{\epsfxsize=2in \centerline{\epsfbox{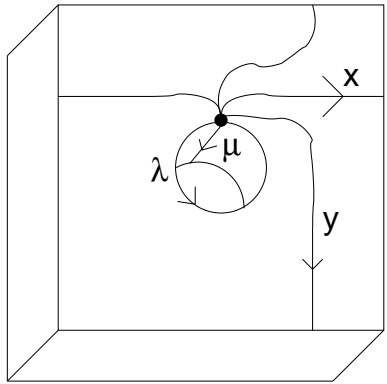}}\hspace{10mm}}
\caption{Notation for the bundle $M$.}
\end{figure}

\section{Homology and boundary slopes of bundles}
\noindent
 Unless otherwise specified, all homology groups in this
 paper will be taken with $\mathbb{Z}$ coefficients.

We begin by recalling some well-known facts.
 If $M$ is a manifold which fibers over $S^1$, with 
 monodromy $f: F \rightarrow F$,
 then there is a corresponding \textit{Wang exact sequence}

$$\begin{CD}
 ... @> >> H_j(M) @> \theta >> H_{j-1}(F) @> f_* - Id >>
 H_{j-1}(F) @> i_* >> H_{j-1}(M) @> >> ...
\end{CD}
$$
where $\theta$ is the map induced by intersection with $F \times \{0\}$.
 There is also a relative version of this sequence, which
 fits with the exact sequence of a pair into the following
 commutative diagram:\\
\\
$\begin{CD}
...  H_j(M, \partial M) @> >> H_{j-1}(F, \partial F) @>  >>
 H_{j-1}(F, \partial F) @>  >> H_{j-1}(M, \partial M) @> >> ...\\
 @VV  V       @VV V   @VV V  @VV V\\
... H_{j-1}(\partial M) @> >> H_{j-2}(\partial F) @>  >>
 H_{j-2}(\partial F) @> >> H_{j-2}(\partial M) @> >> ...\\
 @VV V       @VV V   @VV V  @VV V\\
... H_{j-1}(M) @> >> H_{j-2}(F) @>  >>
 H_{j-2}(F) @>  >> H_{j-2}(M) @> >> ...\\
\end{CD}
$
\\
\\
 In the case where $M$ is a 3-manifold,
 one may use the Wang sequence to compute that
 $Rank H_2(M, \partial M, \mathbb{Q}) = 1 + Rank(fix(f_*),\mathbb{Q})$,
 where $fix(f_*)$ denotes the set of vectors
 in $H_1(F, \partial F)$ which are fixed by $f_*$.

Suppose now that
 $f$ acts trivially on $\partial F$.
 As described in the previous section, for each
 component of $\partial M$, there is a canonical meridian,
 given by the suspension of a point on $\partial F$, and longitude,
 given by intersection with $F$.
 This gives a canonical basis for $H_1(\partial M)$,
 and we have corresponding projections
 $\mu:H_1(\partial M) \rightarrow H_0(\partial F)$,
 which sends each longitude to 0,
 and $\lambda:H_1(\partial M) \rightarrow H_1(\partial F)$,
 which sends each meridian to 0.
 We will fix a preferred component $\ell$ of $\partial F$,
 and let $\pi_{\ell}$ be the projection
 map from $H_1(\partial F)$ to $H_1(\ell)$.
 
 We may define an injective map
 $\eta: H_2(M, \partial M) \rightarrow H_1(F, \partial F) \times \mathbb{Z}$
 by
 $$\eta [R] =([R \cap F \times \{0 \}],
                \pi_{\ell} \lambda([\p R])).$$

 If $[\delta] \in fix(f_*) \subset H_1(F, \partial F)$,
 then $[f \delta - \delta] = 0 \in H_1(F, \partial F)$,
 and so the closed loop $f \delta - \delta$ is homologous
 to some class $[x]$ in
 $Image(i_*:H_1(\partial F) \rightarrow H_1(F)))$.
 Note that $[x]$ is unique, up to adding copies of $[\partial F]$.
 Thus we may
 define $\psi: fix(f_*) \times \mathbb{Z} \rightarrow H_1(\partial F)$
 uniquely, by requiring
 \begin{eqnarray*} 
i_* \psi ([\delta],m) &=& [f\delta - \delta], \textrm{ and}\\
\pi_{\ell} \psi ([\delta],m) &=& m[\ell].
\end {eqnarray*}
 Finally, let $\phi:H_2(M, \partial M) \rightarrow H_1(\partial M)$
 come from the sequence of the pair.

 We have the following diagrams:
$$\begin{CD}
 H_2(M, \partial M) @> \eta  >> fix(f_*) \times \mathbb{Z}\\
 @VV \phi V                  @VV \psi V\\
 H_1(\partial M) @> \lambda >> H_1(\partial F)
\end{CD}
$$
and
$$\begin{CD}
 H_2(M, \partial M) @> >> H_1(F, \partial F)\\
 @VV V                  @VV  V\\
 H_1(\partial M) @> \mu >> H_0(\partial F),
\end{CD}
$$
where in the second diagram, all unlabeled maps come
 from the exact sequence of the pair, and the Wang sequence.
\begin{lem} \label{commutes}
The above diagrams commute.
\end{lem}

\begin{proof}
 For the second diagram,
 note that the map $\mu$ agrees with
 the map $\theta$ from the Wang sequence fpr $\p M$.
 Thus the second diagram fits into the larger diagram
 given at the beginning of the section, which commutes.

 For the first diagram, suppose we are given a class
 $[R] \in H_2(M, \partial M)$.
 A 2-chain homologous to $R$ may  be constructed as follows.
 We let $\delta = R \cap (F \times \{0 \})$,
 and let $\eta[R] = ([\delta], m)$.
 By definition of the map $\psi$,
 we have $[\delta - f \delta] = - i_* \psi ([\delta],m)$.
 Thus there is a map $g:X \rightarrow F$,
 where $X$ is an orientable
 surface with $\partial X = \partial_0 X \amalg \partial_1 X$,
  such that
 $g|_{\partial_0 X}$ is the immersed curve  $\delta - f \delta$ and
 such that $g \partial_1 X \subset \partial F$, with
 $i_* [g \partial_1 X] = i_* \psi ([\delta],m)$;
 since $i_*[\partial F] = 0$, we may also assume that
 $g \partial_1 X \cap \ell = \emptyset$.
 Let $\partial_{00} X \subset \partial_0 X$ be the union of arcs which map to
 $\delta$, and let $\partial_{01} X \subset \partial_0 X$ be the union of arcs
 which map to $-f \delta$.
 Let $\sigma$ be a properly embedded
 collection of separating arcs in $X$, with
 $\partial \sigma = \partial_{00} X \cap \partial_{01} X$.
 Let $h: X- \sigma \rightarrow \{0, 1 \}$ be a continuous map,
 such that $h (\partial_{00} X) = 0$ and $h (\partial_{01} X) = 1$.

 Let $Y$ be an orientable surface
 obtained from $\overline{X - \sigma} \amalg (\sigma \times [0,1])$,
 by identifying
 $\partial_{00} X$ and $\partial_{01} X$ according to the map $f$,
 and identifying $\sigma \times \p [0,1]$ with
 $\overline{X - \sigma} - (X - \s)$ in the obvious way.
 
 Then we may construct a
 map $j: Y \rightarrow M$, 
 by the rule  $j(x) = (g x, h x)$, if $x \in  X -\s$,
 and  $j(x,t) = (g x, t)$ if $(x,t) \in \sigma \times [0,1]$.
 By construction, $\eta([jY]) = ([\d],0)$, and so
 $\eta ([j Y] +m[F]) = ([\delta], m) = \eta [R]$.
 Thus $[j Y] + m[F] = [R]$.
 Also, by inspection,
 $\lambda \phi [j Y] = [g \partial_1 X] =  \psi ([\delta], 0)$,
 and so
\begin{eqnarray*}
 \lambda \phi [R] &=& \lambda \phi ([jY] + m[F])\\
 &=& \psi ([\delta],0)  + m [\partial F]\\
 &=& \psi([\delta],m)\\
 &=& \psi \eta ([j Y] +m[F])\\
 &=& \psi \eta[R].
\end{eqnarray*}
\end{proof}

\begin{cor}
 Suppose $R$ is a properly embedded, orientable, non-separating surface in $M$,
 and let $\delta = R \cap (F \times \{0\})$.
 Then the algebraic boundary slopes of $R$ satisfy:
\begin{eqnarray*} 
\mu([\partial R]) &=& [\delta \cap \partial F] \in H_0(\partial F),\\
 \lambda([\partial R]) &=& [f \delta - \delta] + k[\p F] \in H_1(\p F)\,\,\,
 \textrm{ for some integer } k. 
\end{eqnarray*}

\end{cor}

The following corollary employs notation introduced in Section 2.

\begin{cor} \label{slopes}
Suppose there is an arc $\delta$ properly embedded in $F$,
 with $[\delta \cap \partial F] = \sum a_i [p_i] \in H_0(\partial F)$,
 and $[f \delta - \delta] = \sum b_i [\lambda_i] \in H_1( F)$.
 Then there is a non-separating, orientable surface
 $R$ properly embedded in $M$,
 with algebraic boundary slopes $((a_1, b_1), ..., (a_k, b_k))$.
\end{cor}

\begin{proof}
By the Wang exact sequence, there is a class
 $[R] \in H_2(M, \partial M)$, such that
 $[R \cap F \times \{0 \}] = \theta [R] = [\delta] \in H_1(F, \partial F)$.
  Then, by the previous corollary,
 we have $\mu([\partial R]) = [\delta \cap \partial F] \in H_0(\partial F)$,
 and $\lambda([\partial R]) = [f \delta - \delta] + k[\p F]
 \in H_1(\partial F)$.
 Adding a multiple of $[F]$ to $[R]$, we get
 $\l[\p R] = [f \d - \d]$, and the corollary follows.
\end{proof}

\section{Proof of Theorem \ref{main}}

We shall make use of the following,
 which can be proved by straightforward applications
 of the methods of \cite{CL}.

\begin{thm} \label{bc}
 For any compact,
 orientable surface, $S$,  there is a positive integer $n = n(S)$,
 depending only on the
 topological type of $S$, such that the following is true.
 Let $M$ be any compact, orientable, irreducible,
 atoroidal 3-manifold, with two torus boundary components,
 containing a properly embedded, orientable, incompressible,
 non-separating surface homeomorphic to $S$,
 which is not a fiber in a fibration of $M$,
 with algebraic boundary slopes $(\d_1, \d_2)$,
 where $\delta_i \neq 0$.
 Then, if $|I(\alpha_1, \delta_1)| = |I(\alpha_2, \delta_2)| > n(S)$,
 the manifold obtained by Dehn filling $M$ along
 $\alpha_1$ and $\alpha_2$ is virtually Haken.
\end{thm}
 
\textit{Remark:} By choosing the function $n(S)$ appropriately,
 we may assume that, whenever $S^{\prime}$ is obtained
 by compressing $S$, we have $n(S^{\prime}) < n(S)$ .\\
\\
We are now ready to prove the main theorem.
\begin{proof} (Of Theorem \ref{main})

 It is well-known that, if $f, g \in \mathcal{M}_1^1$,
 and  $f_* = g_*$, then $f g^{-1}$ is isotopic to a power of a
 Dehn twist along a peripheral curve in $F$.
 Therefore, if $f_* = g_*$, then $M_f$ is bundle equivalent to $M_g$,
 and so we may reduce to the case where
 $f \in J_i = <D_y^i, D_x>$.\\
\\
 Case a: $f \in J_3$.\\

 In this case, $f$ can be written as a word $W$ in
 $D_x$ and $D_y^3$.
 Let $m = m_3$ denote the exponent sum of $D_x$ in $W$.
 Let $\beta = \beta_3^1$ be the slope  $(3, m)/gcd(3,m)$ on $\partial M$.

 Let $\theta: \widehat{F} \ra F$
 be the 3-fold cyclic cover of $F$ dual to the curve $x$
 (see Fig. 1).
 Since $f \in H_3$, then $f$ lifts to an automorphism
 $\widehat{f} : \widehat{F} \rightarrow \widehat{F}$
 which acts trivially on $\partial \widehat{F}$. 
 Let $\pi:\widehat{M} \rightarrow M$ be the 3-fold cover of $M$ induced
 by the lift $\widehat{f}: \widehat{F} \rightarrow \widehat{F}$.

 If $g$ is a map between two surfaces with boundary,
 we let $g_{\sharp}$ be the induced map on
 $H_1$ rel. boundary. 
   Let $\lambda_1, \lambda_2, \lambda_3$ be the components
 of $\partial \widehat{F}$
 (with pre-image orientations induced from $\lambda$),
 and let $p_i \in \lambda_i$ be the pre-images
 of $p \in \partial F$.
 Let $\delta_i$ be an arc connecting $\lambda_i$ and
 $\lambda_{i+1}$, as pictured in Figure 2a.
 Let $\delta = \delta_1 - \delta_2$.
 Then it is easy to see that $[\delta]$ is
 a non-zero class in $H_1(\widehat{F}, \partial \widehat{F})$
 which is fixed by the lifts
 of $D_{x}$ and $(D_y)^3$, and therefore by
 $\widehat{f}_{\sharp}$.
 Moreover, $[f \delta - \delta] = -m[\lambda_2] \in H_1(F)$.

 \begin{figure}[!htbp]
{\epsfxsize=4in \centerline{\epsfbox{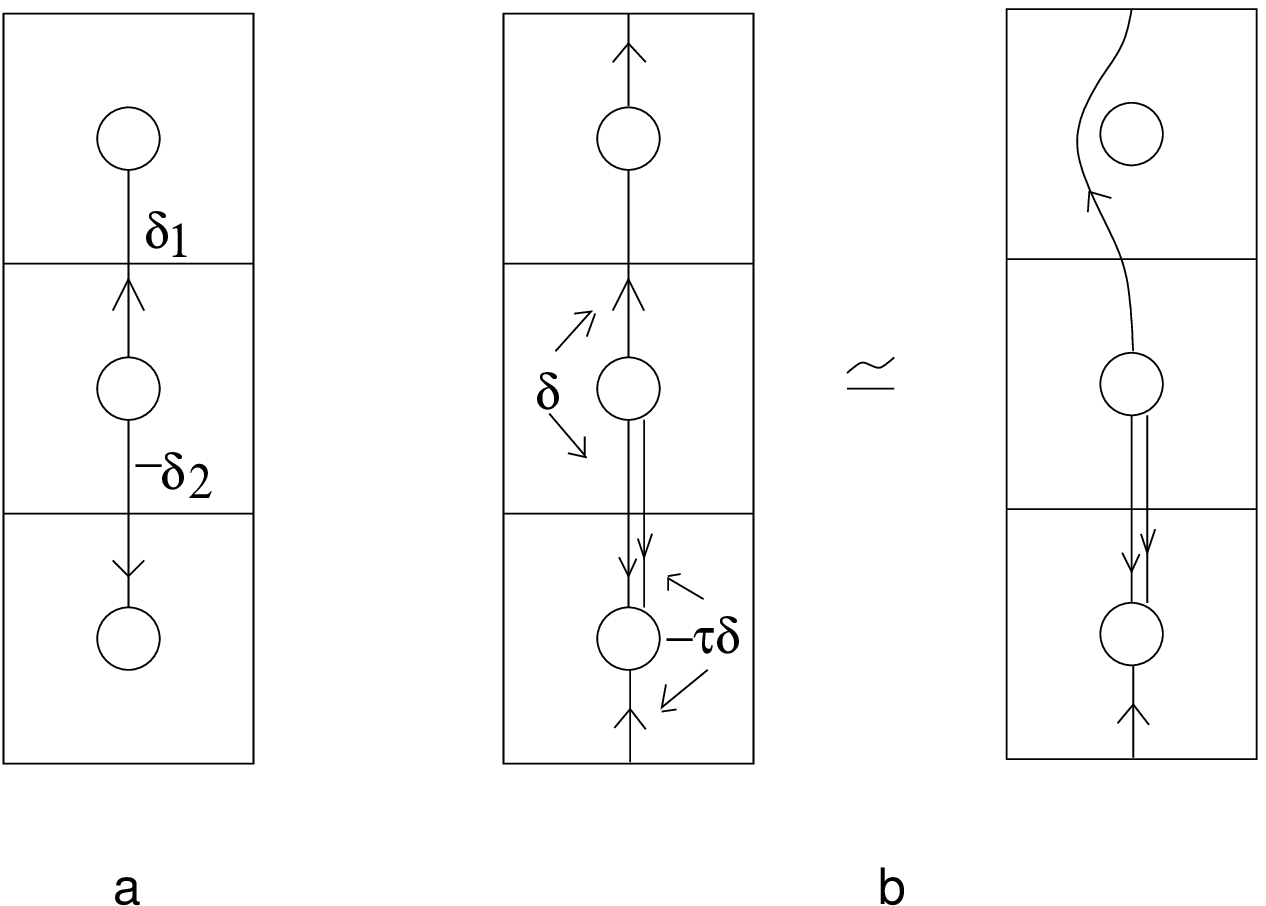}}\hspace{10mm}}
\caption{a. A class in $ker(p_*)$ \,\,
 b. A fixed class in $H_1(\widehat{F}, \partial \widetilde{F})$,
 corresponding to the surface $S$.}
\end{figure}

 Let $\tau$ be the covering transformation of $\widehat{F}$
 such that $\tau (\delta_1) = \delta_2$.

\begin{lem} 
There is a non-separating surface $R \subset M$ 
 whose boundary slopes are given by $((0,0), (-3, -m), (3,m))$.
 \end{lem}

\begin{proof}
 We have
 $[f (\delta - \tau \delta) - (\d - \t \d)] =
  m [\l_3] - m[\lambda_2] \in H_1(F)$,
 and $[(\delta - \t \d) \cap \partial F] = 0[p_1] - 3[p_2] + 3[p_3]
 \in H_0(\partial F)$.
 The lemma now follows from Corollary \ref{slopes}.
\end{proof}

Let $T_i$ be the component of $\p \wh M$ containing $p_i$.
 By attaching annuli to $R$, we may assume that $R$ is disjoint
 from $T_1$.
Let $\alpha$ be a slope on $\partial M$,
 let $\widehat{\alpha}_i$ be the lift of $\alpha$ to $T_i$,
 and let $S$ be an incompressible surface in
  $\widehat{M}(\widehat{\alpha}_1)$ obtained by compressing $R$.
 Then $S$ has boundary slopes which project
 to curves of slope $\beta$ in $\partial M$.

Let $N$ be the infinite
 cyclic cover of $\wh M(\widehat{\alpha}_1)$
 which is dual to $S$, and let $\Delta(s)$ be
 the Alexander polynomial associated to this infinite
 cyclic cover. 

 Let $\th: \wt F \ra \wh F$ be the infinite cyclic cover
  dual to $\d - \t \d$. 
 We let $\wh x_i = y^{i-1} x y^{-(i-1)} \in \wh F$, and
 let $\wh y = y^3 \in \wh F$.
 There is a $\mathbb{Z}[s^{\pm 1}]$-module decomposition
 $H_1(\wt F, \th^{-1} p_1) \cong H_1(\wt F) + <[\wh x_1]>$.

 The monodromy $\wh f$ acts on the
 module $H_1(\wt F, \theta^{-1} p_1 )$
 by a matrix $f_*$, and since $\wh f$ fixes $p_1$,
 then $Im(Id - f_*) \subset H_1(\wt F)$.
 The map $Id - f_*$ has a kernel containing $<[\wh x_1 \wh x_2^{-1}]>$,
 and so there is an induced map
 $\overline{Id - f_*}:  H_1(\wt F, \th^{-1} p_1)/<[\wh x_1 \wh x_2^{-1}]>
                            \ra H_1(\wt F) \subset H_1(\wt F, \th^{-1} p_1)$.
 
 The group $H_1(\wt F, \th^{-1} p_1)$ is a free
 $\mathbb{Z}[s^{\pm 1}]$-module, with basis
 $\mathcal{B} =
  ([\wh y], [\wh x_1 \wh x_2^{-1}], [\wh x_1^2 \wh x_3], [\wh x_1])$.
  We let
 $[f_*]$ be the matrix representative for $f_*$ with respect
 to $\mathcal{B}$.  Then $\overline{Id - f_*}$ can
 be represented by a matrix $[\overline{Id - f_*}]$,
 which is obtained by deleting the second column and fourth row
 of $Id - [f_*]$.

 \begin{lem} \label{alex1}
a. The matrix $[f_*]$ is given by
 matrix $W([D_{x*}], [D_{y*}^3])$, where
\begin{eqnarray*}
[D_{x*}] =
 \begin{pmatrix}
s & 0 & 0 & 0\\
s^{-1} & 1 & s-s^{-1} & 0\\
-s^{-1} & 0 & s^{-1} & 0\\
0 & 0 & 0 & 1
\end{pmatrix},
\end{eqnarray*}
 and

\begin{eqnarray*}
[D_{y*}^3] =
 \begin{pmatrix}
1 & 0 & 1+s+s^2 & 1 \\
0 & 1 & 0 & 0 \\
0 & 0 & 1 & 0\\
0 & 0 & 0 & 1
\end{pmatrix},
\end{eqnarray*}

b. The Alexander polynomial $\D(s)$ is the determinant of the following
 matrix:
\begin{eqnarray*}
B = \begin{pmatrix}
0 &  & &  \\
0 &  & [\overline{Id-f_*}]^T & \\
1-s &  &  & \\
p & 0 & -qs^{-1} & 0
\end{pmatrix}.
\end{eqnarray*}
\end{lem}

Using Lemma \ref{alex1} and a computer,
 the Alexander polynomial $\D(s)$ can be easily computed.

\begin{proof}

 Part a is a computation, which can be done with Fox derivatives.
 We leave this to the reader.

 For part b,
 let $\widetilde{M}$ be the $\mathbb{Z}$-cover of $\widehat{M}$
 which is dual to $R$.
 We begin by computing a presentation for $H_1(\widetilde{M})$ as
 a $\mathbb{Z}[s^{\pm 1}]$- module.
 
 We choose the point $p_1$ as a basepoint for $\pi_1 \wh M$.
 We have $\pi_1 \wh F \subset \pi_1 \wh M$, and we let
 $\wh t = t \subset \pi_1 \wh M$.
 We have the following presentation for $\pi_1(\widehat{M})$:

\begin{eqnarray*}
\pi_1 \wh M = < \wh y, \wh x_1, \wh x_1 \wh x_2^{-1},
 \wh x_1^2 \wh x_3, \wh t|R_1, R_2, R_3, R_4>,\\
R_1= \wh t \wh y \wh t^{-1} \wh y^{-1} (\wh y \wh f_* \wh y^{-1}),\\
R_2 = \wh t \wh x_1 \wh t^{-1} \wh x_1^{-1} (\wh x_1 \wh f_* \wh x_1^{-1})\\
R_3=  \wh t \wh x_1 \wh x_2^{-1} \wh t^{-1} (\wh x_1 \wh x_2^{-1})^{-1}
 (\wh x_1 \wh x_2^{-1} \wh f_* (\wh x_1 \wh x_2^{-1})^{-1}),\\
R_4= \wh t \wh x_1^2 \wh x_3 \wh t^{-1} (\wh x_1^2 \wh x_3)^{-1}
 (\wh x_1^2 \wh x_3 \wh f_* (\wh x_1^2 \wh x_3)^{-1}).
\end{eqnarray*}

 For any element $w \in \pi_1 \wt M \subset \pi_1 M$,
 it will be convenient to let $[w]$ denote the the image
 of $w$ in $H_1(\wt M)$, and to
 let $s = \wh x_1$.  Then, as a $\mathbb{Z}[s^{\pm 1}]$-module,
 $H_1(\wt M)$ has an ordered generating set 
$\mathcal{B}^{\prime} =
 ([\wh t], [\wh y], [\wh x_1 \wh x_2^{-1}], [\wh x_1^2 \wh x_3])$.

 The relators of $H_1(\wt M)$ (as a $\mathbb{Z}[s^{\pm 1}]$-module)
 may be obtained from the relators $R_1, ..., R_4$ of $\pi_1 \wt M$.
 We have:
\begin{eqnarray*}
R_1 &=&\wh t \wh y \wh t^{-1} \wh y^{-1} (\wh y \wh f_* \wh y^{-1}),\\
 & \Rightarrow& [\wh y \wh f_* \wh y^{-1}] = 0\\
R_2 &=&  \wh t \wh x_1 \wh t^{-1} \wh x_1^{-1} (\wh x_1 \wh f_* \wh x_1^{-1})\\
 & \ra& (1-s)[\wh t] + [\wh x_1 \wh f_* \wh x_1^{-1}] = 0\\
R_3 &=&  \wh t (\wh x_1 \wh x_2^{-1}) \wh t^{-1} (\wh x_1 \wh x_2^{-1})^{-1}\\
 &\Rightarrow& 0 = 0\\
R_4 &=&  \wh t (\wh x_1^2 \wh x_3) \wh t^{-1} (\wh x_1^2 \wh x_3)^{-1}
 ((\wh x_1^2 \wh x_3) \wh f_* (\wh x_1^2 \wh x_3)^{-1})\\
 &\Rightarrow& [(\wh x_1^2 \wh x_3) \wh f_* (\wh x_1^2 \wh x_3)^{-1}] = 0\\
\end{eqnarray*}

Thus a presentation matrix for $H_1(\wt M)$,
 in terms of $\mathcal{B}^{\prime}$,
  is given by:
\begin{eqnarray*}
A = \begin{pmatrix}
0 &  & &  \\
0 &  & \overline{Id-f_*}^T & \\
1-s &  &  & 
\end{pmatrix}.
\end{eqnarray*}
The presentation matrix for $H_1(N)$
 is obtained from $A$ by adding a single
 relator of the form $[t^p (x^{-1} y x y^{-1})^q] = 0$,
 where $p, q \in \mathbb{Z}$.
 This yields the relator $p[\wh t] - qs^{-1} [x_1 x_2^{-1}] = 0$.
 Thus the presentation matrix for  $H_1(N)$ is:
\begin{eqnarray*}
\begin{pmatrix}
0 &  & &  \\
0 &  & \overline{Id - f_*}^T & \\
1-s &  &  & \\
p & 0 & -qs^{-1} & 0
\end{pmatrix},
\end{eqnarray*}
and so $\D(s)$ is the determinant of this matrix.
\end{proof}

 \begin{cor} \label{nonfiber1}
  There is a slope $\beta_3^2$ such that, if $I(\alpha, \beta_3^2) > 1$,
 then $S$ is not a fiber
 in a fibration of $\widehat{M}(\widehat{\alpha}_1)$.
 \end{cor}

\begin{proof}
 It is well-known that, if an infinite cyclic cover
 of a compact 3-manifold is dual to a fiber in a fibration,
 then the corresponding Alexander polynomial is monic.
 By Lemma \ref{alex1}, $\Delta(s) = p Det M_1 + q Det M_2$
 for some matrices $M_1$ and $M_2$ with entries in
 $\mathbb{Z}[s^{\pm 1}]$.
 Thus the leading term of $\Delta(s)$ is $p n_1  + q n_2$,
 for some integers $n_1, n_2$. We let
 $\b_3^2 = (-n_2, n_1)/gcd(n_2,n_1)$.
 Then if $|I(\alpha, \b_3^2)| \neq 1$, $\D(s)$ is non-monic. 
\end{proof}

  Let $n$ be a positive integer.
 By Thurston's hyperbolic Dehn surgery theorem,
 we may assume that $n$ is chosen large enough so that 
 $\widehat{M}(\widehat{\alpha}_1)$ is hyperbolic whenever
 $I(\alpha, \beta_3^1) > n$.   We also assume that
 $n$ is larger than the integer $n(R)$ given by Theorem \ref{bc},
 and hence also bigger than $n(S)$.

 Suppose $|I(\alpha, \beta_3^1) | = k > n$,
 and that $|I(\a, \b_3^2)| > 1$.
 Since $k > n$, then $\widehat{M}(\widehat{\alpha}_1)$
 is hyperbolic, and since $|I(\alpha, \b_3^2)| > 1$
 then, by Corollary \ref{nonfiber1}, $S$ is not a fiber in a fibration of
 $\widehat{M}(\widehat{\alpha}_1)$.
 Let $\widehat{\alpha}_2, \widehat{\alpha}_3$
 be the lifts of $\alpha$ to the 
 components of $\widehat{M}(\widehat{\alpha}_1)$.
 Then $|I(\widehat{\alpha}_2, \partial S)|
  = |I(\widehat{\alpha}_3, \partial S)| = k
 > n \geq n(S)$,
 so by Theorem \ref{bc}, the Dehn filling of
 $\widehat{M}(\widehat{\alpha}_1)$
 along $\widehat{\alpha}_2$ and $\widehat{\alpha}_3$
 is virtually Haken.
 Since this manifold covers $M(\alpha)$, then $M(\alpha)$ is virtually Haken.\\
\\
\\
\textit{Case b:} $f \in J_4$.\\
 The argument is similar to the argument for Case a.
 Let $\beta = \beta_4 = (2, m)/gcd(2, m)$, where $m = m_4$ is the exponent
 sum of $D_x$ in the word $W(D_x, D_y^4)$ which represents $f$.
 We let $\widehat{F}$ be the 4-fold cyclic cover dual to $x$,
 and let $\widehat{M}$ be the corresponding cover of $M$.
 We fix a basepoint $p_1 \in \p \wh M$,
 and let $p_i$ be the translate of $p_1$ by the covering
 translation corresponding to $y^{i-1}$.  We let $T_i$ be
 the component of $\p \wh M$ containing $p_i$.

  We define arcs $\delta_1$ and $\delta_2$ in $\widehat{F}$
 as pictured in Figure 3, and let $\delta = \delta_1 - \delta_2$.
  In this case, the
 non-separating surface $R$ corresponding to $[\delta] - [\tau^2 \delta]$
 has boundary slopes
 $((0,0), (-2, -m), (0,0), (2,m)) \neq ((0,0),(0,0),(0,0),(0,0))$.
 By compressing $R$, 
 we obtain a properly embedded, orientable,
 non-separating, incompressible  surface $S$ in
 $\wh M(\wh \a_1, \wh \a_3)$,
 whose boundary curves
 all project to curves of slope $\b$ on $\partial M$.

\begin{figure}[!htbp]
{\epsfxsize=4in \centerline{\epsfbox{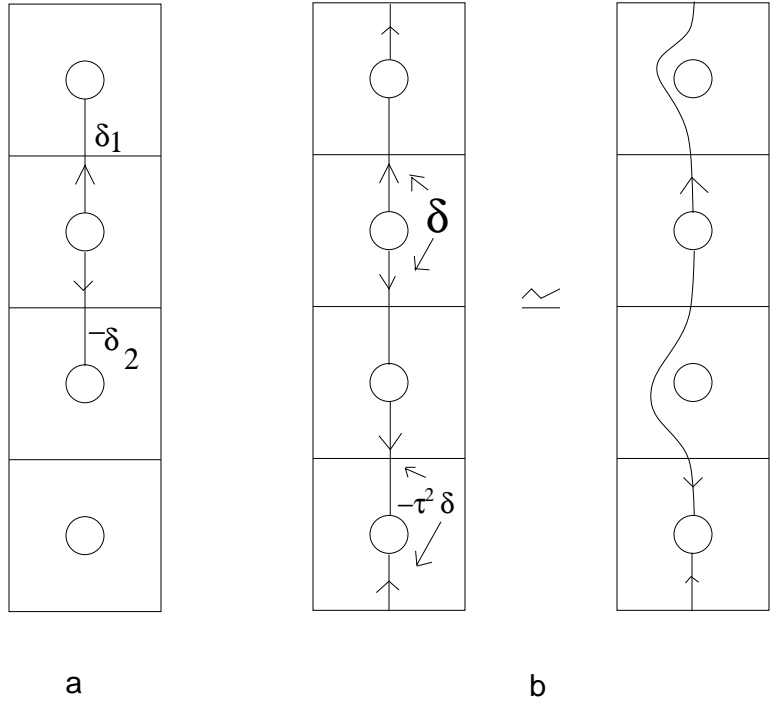}}\hspace{10mm}}
\caption{a. A class in $ker(p_*)$ \,\,\,\,\,\,
 b. A fixed class in $H_1(\widehat{F}, \partial \widetilde{F})$, corresponding to the surface $S$.}
\end{figure}

Let $N$ be the infinite cyclic cover of $\wh M(\wh \a_1, \wh \a_3)$
 dual to $S$,
 and let $\D(s)$ be the corresponding
 Alexander polynomial.
 Let $\th: \wt F \ra \wh F$ be the infinite cyclic cover of $\wh F$ dual
 to $\d - \t^2 \d$.
 We let $\wh x_i = y^{i-1} x y^{-(i-1)} \in \wh F$, and
 let $\wh y = y^4 \in \wh F$.

 The automorphism $\wh f$ induces an automorphism
 $f_*$ of  $H_1(\wt F, \th^{-1} p_1 \cup p_3)$.
 There is an induced map
 $$\overline{Id - f_*} : H_1(\wt F, \th^{-1}(p_1 \cup p_3))/
                   <[\wh x_1 \wh x_2^{-1}], [\wh x_3 \wh x_4^{-1}]>
               \ra H_1(\wt F) \subset H_1(\wt F, \th^{-1} (p_1 \cup p_3)).$$
 The $\mathbb{Z}[s^{\pm 1}]$-module $H_1(\wt F, \th^{-1} (p_1 \cup p_3))$
 is free, with basis
  $$\mathcal{B} = ([\wh y^4], [\wh x_1 \wh x_3], [\wh x_1],
 [\wh y^2], [\wh x_1 \wh x_2^{-1}], [\wh x_3 \wh x_4^{-1}]).$$
 We let $[f_*]$ be the matrix representing $f_*$, in terms
 of the basis $\mathcal{B}$.
 Then $\overline{Id - f_*}$ is represented by a matrix
 $[\overline{Id - f_*}]$ obtained from $Id - [f_*]$
 by deleting the 5th and 6th columns, and the 3rd and 4th rows.

\newpage
\begin{lem} \label{alex2}
a. The matrix $[f_*]$ is given by $W([D_{x*}], [D_{y*}^4]$, where
\begin{eqnarray*}
[D_{x*}] =
 \begin{pmatrix}
s & 0 & 0 & 0 & 0 & 0 \\
-(1+s^{-1}) & s^{-1} & 0 & -s^{-1} & 0 & 0\\
0 & 0 & 1 & 0 & 0 & 0\\
0 & 0 & 0 & 1 & 0 & 0\\
s^{-1} & 1-s^{-1} & 0 & s^{-1} & 1 & 0\\
0 & 0 & 0 & 0 & 0 & 1
\end{pmatrix},
\end{eqnarray*}
 and

\begin{eqnarray*}
[D_{y*}^4] =
 \begin{pmatrix}
1 & 1+s & 1 & 0 & 0 & 0\\
0 & 1& 0 & 0 & 0 & 0\\
0  & 0 & 1 & 0 & 0 & 0\\
0 & 0 & 0 & 1 & 0 & 0\\
0 & 0 & 0 & 0 & 1 & 0\\
0 & 0 & 0 & 0 & 0 & 1
\end{pmatrix},
\end{eqnarray*}

b. The Alexander polynomial $\D(s)$ is the determinant of the following
 matrix:

\begin{eqnarray*}
B = \begin{pmatrix}
0 & 0 &0 &0 &0 \\
0 & 0 &0 &0 &0\\
1-s &0&0 &0 &0 \\
p & 0 & 0 & 0 & -qs\\
p & 0 & 0 & -qs^{-1} & 0
\end{pmatrix}
 + \begin{pmatrix}
0 &  & & & \\
0 &  & \overline{Id - f_*}^T & &\\
0 &  &  && \\
0 &  &  &  & \\
0 & 0 & 0 & 0 & 0
\end{pmatrix}
\cdot
  \begin{pmatrix}
1 & 0 &0 &0 &0 \\
0 & 1 &0  &0 &0\\
0 & 0 & 1 &0&0 \\
0 & 0 &0  & p &0 \\
0 & 0 & 0 & 0 & 1
\end{pmatrix},
\end{eqnarray*}
and $\D(s)$ is divisible by $q$.
\end{lem}

\begin{proof}
 The proof of part a is an elementary application of Fox calculus.
 
For part b,
 let $\wt M$ be the infinite cyclic cover of $\wt M$ dual to $S$.
 We choose an ordered generating set
 $\mathcal{B}^{\prime} = (\wh t, \wh y^4, \wh x_1 \wh x_3,
 \wh x_1 \wh x_2^{-1}, \wh x_3 \wh x_4^{-1})$ for
 the $\mathbb{Z}[s^{\pm 1}]$-module $H_1(\wt M)$.

 The relators for $\pi_1 \wh M$ give the following relations
 for $H_1(\wt M)$:

\begin{eqnarray*}
R_1 &=&\wh t \wh y \wh t^{-1} \wh y^{-1} (\wh y \wh f_* \wh y^{-1}),\\
 & \Rightarrow& [\wh y \wh f_* \wh y^{-1}] = 0\\
R_2 &=&  \wh t \wh x_1 \wh x_2^{-1} \wh t^{-1} (\wh x_1 \wh x_2^{-1})^{-1}
 (\wh x_1 \wh x_2^{-1} \wh f_*  (\wh x_1 \wh x_2^{-1})^{-1})\\
 &\Rightarrow& 0 = 0\\
R_3 &=&  \wh t \wh x_1 \wh x_3 \wh t^{-1} (\wh x_1 \wh x_3)
 (\wh x_1 \wh x_3 \wh f_*  (\wh x_1 \wh x_3)^{-1})\\
 &\Rightarrow& [\wh x_1 \wh x_3 \wh f_* (\wh x_1 \wh x_3)^{-1}] = 0\\
R_4 &=&  \wh t \wh x_1 \wh t^{-1} \wh x_1^{-1} (\wh x_1 \wh f_* \wh x_1^{-1})\\
 & \ra& (1-s)[\wh t] + [\wh x_1 \wh f_* \wh x_1^{-1}] = 0\\
R_5 &=&  \wh t \wh x_3 \wh x_4^{-1} \wh t^{-1} (\wh x_3 \wh x_4^{-1})
 (\wh x_3 \wh x_4^{-1} \wh f_*  (\wh x_3 \wh x_4^{-1})^{-1})\\
 &\Rightarrow&  0 = 0\\
\end{eqnarray*}

The module $H_1(N)$
 has two additional relators,
 $R_6: [t^p (x^{-1} y x y^{-1})^q] = 0$, and
 $R_7: [y^2t^py^-2 y^2(x^{-1} y x y^{-1})^qy^{-2}] = 0$.
 The relator $R_6$ can be written as:
$p [\wh t] -qs^{-1} [\wh x_1 \wh x_2^{-1}] = 0$.
 For  $R_7$, we have:
\begin{eqnarray*}
0 &= & p[y^2  t y^{-2}] + q[y^2(x^{-1} y x y^{-1})y^{-2}]\\
 &=& p[(y^2 t y^{-2} t^{-1}) t] -qs [x_1 x_2^{-1}]\\
 &=& p[\wh t] + p[y^2 f_* y^{-2}] - qs [x_1 x_2^{-1}]\\
 &=& p[\wh t] + -p[\overline{Id -f_*} \,\, y^2] - qs [x_1 x_2^{-1}].
\end{eqnarray*}
 
 Thus the presentation matrix for  $H_1(N)$, with respect to $\mathcal{B}$, is:
\begin{eqnarray*}
\begin{pmatrix}
0 & 0 &0 &0 &0 \\
0 & 0 &0 &0 &0\\
1-s &0&0 &0 &0 \\
p & 0 & 0 & 0 & -qs\\
p & 0 & 0 & -qs^{-1} & 0
\end{pmatrix}
 + \begin{pmatrix}
0 &  & & & \\
0 &  & \overline{Id - f_*} & &\\
0 &  &  && \\
0 &  &  &  & \\
0 & 0 & 0 & 0 & 0
\end{pmatrix}
\cdot
  \begin{pmatrix}
1 & 0 &0 &0 &0 \\
0 & 1 &0  &0 &0\\
0 & 0 & 1 &0&0 \\
0 & 0 &0  & p &0 \\
0 & 0 & 0 & 0 & 1
\end{pmatrix}
\end{eqnarray*}

We then observe that the last column of $\overline{Id - f_*}$
 is all 0's. Therfore the last column of
 the presentation matrix for $H_1(N)$
 has a single non-zero entry $-q s$.  Thus $\D(s)$ is divisible by $q$.
\end{proof}

 \begin{cor} \label{nonfiber2}
 Let $\b_4^2 = (1,0)$.
 If $I(\a, \b_4^2) >1$, then $S$ is not a fiber
 in a fibration of $\widehat{M}(\widehat{\alpha}_1, \wh \a_3)$.
 \end{cor}

\begin{proof}
By Lemma \ref{alex2}, $\D(s)$ is divisible by $q$. 
 If $|q| = |I(\a, (1,0))|>1$,
 then $\D(s)$ is non-monic.
\end{proof}

 Now an application of Theorem \ref{bc} shows that
 $\widehat{M}(\widehat{\alpha}_1, \widehat{\alpha}_2,
 \widehat{\alpha}_3, \widehat{\alpha}_4)$
 has a Haken cyclic cover, provided that $|I(\alpha, \beta_4^2)| > 1$
 and $|I(\a, \b_4^1)|$ is large enough.
\end{proof}

\section{Computations: framings}

In the proof of Theorem \ref{main},
 we used the fact
 that, if $f,g \in \mathcal{M}_1^1$, and $f_* = g_* \in SL_2(\mathbb{Z})$,
 then $M_f$ is bundle-equivalent to $M_g$.
 In this section, we shall show how to compute
 the effect of this equivalence on the framings.

 Let $D_{\lambda}$ be a Dehn twist about a peripheral curve in $F$.
 Suppose $f \in \mathcal{M}_1^1$, and that $f_* = Id$.
 Then $f$ is equivalent in $\mathcal{M}_1^1$ to $D_{\lambda}^n$ for
 some integer $n$. We define the \textit{twist} of $f$
 by the formula $t(f) = n$.

\begin{lem} \label{frame}
Suppose that, $f, g \in \mathcal{M}_1^1$, with $f_* = g_*$.
 Then there is a homeomorphism $h:M_f \rightarrow M_g$,
 such that, with respect to the standard
 framings on $\partial M_f$ and $\partial M_g$,
 $h(1,0) = (1, t(f g^{-1}))$ and $h(0,1) = (0,1)$.
\end{lem}

\begin{proof}
 Since $f g^{-1} = D_{\lambda}^{t(f g^{-1})}$,
 then $f g^{-1}$ is isotopic to the identity, by an isotopy
 which twists $t(f g^{-1})$ times around the boundary of $F$.
 Using this isotopy, one may construct a bundle-equivalence
 between $M_f$ and $M_g$, and verify that the effect on the framings
 is as claimed.
\end{proof}

 Thus, given elements $f, g \in \mathcal{M}_1^1$, in terms of $D_x$ and $D_y$,
 with $f_* = g_*$,
 we require a method for computing the twist $t(f g^{-1})$.

\begin{lem} \label{twist}
Let $f, g \in \mathcal{M}_1^1$ be given as words,
 $W_f, W_g$ in $D_x$ and $D_y$, and suppose that $f_* = g_*$.
 Then $t(f g^{-1})$ is equal to 1/12 of the total sum of the
 exponents of $D_x$ and $D_y$ in
 the word $f g^{-1} = W_f(D_x, D_y) W_g^{-1}(D_x, D_y)$.

\end{lem}

\begin{proof}
 We use the following well-known presentation for
 $SL_2(\mathbb{Z})$:
$$ SL_2(\mathbb{Z}) \cong <a, b, \tau
    | \tau = (a b)^3 = (a b a)^2, \tau^2 = [\tau, a] = [\tau, b]= id>,$$
 where the map sending $D_{x*}$ to $a$ and $D_{y*}$ to $b$ is an isomorphism.
 Thus if $f_* = g_*$, then $fg^{-1} = W(D_x, D_y)$ is a product of conjugates
 of the elements $R_1 = (D_x D_y)^3(D_x D_y D_x)^{-2}$,
 $R_2 = (D_x D_y)^3 D_x (D_x D_y)^{-3} D_x^{-1}$,
  $R_3 = (D_x D_y)^3 D_y (D_x D_y)^{-3} D_y^{-1}$,
 and $R_4 = (D_x D_y)^6$.
 By computing the effect of these automorphisms on $\pi_1 F$,
 one may check directly that each one is trivial
 in $\mathcal{M}_1^1$ except for $R_4$,
 and that $t(R_4) = 1$.
 Therefore $t(f g^{-1})$ is equal to the (signed) number of conjugates
 of $R_4$ in $W_f(D_x, D_y) W_g^{-1}(D_x, D_y)$.
 Since the sum of the exponents is 12 on $R_4$, and
 zero on each of $R_1, R_2$ and $R_3$,
 then we see that $t(f g^{-1})$ is simply 1/12 of the the total sum of the
 exponents of $W_f W_g^{-1}$.
 \end{proof}

\section{Computations: subgroups of $SL_2(\mathbb{Z})$}

The following lemma is known.  It can
 be proved directly, using standard combinatorial group theory algorithms
 (as implemented for example on the program GAP).
 We shall give a proof based on the Euclidean algorithm.
 
\begin{lem}
The subgroups $H_3$ and $H_4$ have finite index
 in $\mathcal{M}_1^1$.
\end{lem}

\begin{proof}
We shall prove the statement for $H_4$, the proof for $H_3$
 being entirely analogous.

Given an ordered pair of relatively prime, non-zero integers $(m,n)$,
 we may generate a sequence $(m_0, n_0), ..., (m_k, n_k)$
 recursively, as follows:\\
Let $(m_0, n_0) = (m,n)$.
Suppose $(m_i, n_i)$ has been defined.
 If $|n_i \pm m_i| < |n_i|$, then
 let $(m_{i+1},n_{i+1}) = (m_i,n_i \pm m_i)$;  if $|n_i + m_i| \geq |n_i|$
 and $|n_i - m_i| \geq |n_i|$ and
 $|m_i \pm 4 n_i| < |m_i|$, then let
 $(m_{i+1},n_{i+1}) = (m_i \pm 4n_i, n_i)$;
 if neither of these conditions holds, then terminate
 the sequence at $(m_i,n_i)$.\\
\\
\textit{Claim:} For any pair of relatively prime, non-zero integers
 $(m,n)$, the above rule defines a finite sequence 
 terminating in $(\pm 1,0)$ or $(0,\pm 1)$.\\
\\
\textit{Proof of claim:}
Suppose that $(m_i, n_i)$ has been defined, that neither
 $m_i$ nor $n_i$ is zero, and that $|n_i + m_i|$ and
 $|n_i - m_i|$ are both as big as $|n_i|$.
  Then $|m_i| \geq 2|n_i|$.  Since $m$ and $n$ are relatively
 prime, it follows that $m_i$ and $n_i$ are relatively prime, and thus
 we have strict inequality $|m_i| > 2|n_i|$,
  and so $|m_i \pm 4n_i| < |m_i|$.
 So the sequence continues until
 we reach $(j,0)$ or $(0,j)$.
 In this case, both $m$ and $n$ are divisible
 by $j$, so we have $j= \pm 1$. This proves the claim.\\
\\
  Now, $SL_2(\mathbb{Z})$ acts on the hyperbolic plane,
 and there is an induced action on the circle
 at infinity, which is identified with $\mathbb{R} \cup \{ \infty \}$.
 
Identifying the rational number
 $m/n$ with the vector $ \begin{pmatrix}
m \\
n
\end{pmatrix}$, and $\infty$ with $ \begin{pmatrix}
1 \\
0
\end{pmatrix}$, the action of $<D_{x*}, D_{y*}^4>$ on
 $\mathbb{Q} \cup {\infty}$ is given by:
 
$\begin{pmatrix}
1 & \pm 1\\
0 & 1
\end{pmatrix}
\begin{pmatrix}
m \\
n
\end{pmatrix}
 =\begin{pmatrix}
m \pm n\\
n
\end{pmatrix}
$

 and
 $\begin{pmatrix}
1 & 0\\
\pm 4 & 1
\end{pmatrix}
\begin{pmatrix}
m \\
n
\end{pmatrix}
 =
\begin{pmatrix}
 m\\
\pm 4m +n
\end{pmatrix}
.$
Therefore, by the claim, the orbit
 $<D_{x*}, D_{y*}^4> \{0, \infty \}$
 is dense in $\mathbb{R} \cup \{ \infty \}$.
  Thus the domain of discontinuity of
  $<D_{x*}, D_{y*}^4>$  is empty, and
 so  $<D_{x*}, D_{y*}^4>$ is a finite-index subgroup
 of $SL_2(\mathbb{Z})$.  Thus $H_4$ has finite index in $\mathcal{M}_1^1$.
\end{proof}

In fact, it can be checked (for example on  GAP),
 that $H_3$ has index 8, and $H_4$ has index 12.\\
\\
\textit{Computation of exponent sums:}

 Given an element $g \in \mathcal{M}_1^1$, let $m_i$ be the smallest
 positive integer such that $g^{m_i} \in H_i$.
 Thus there is word $W$ (not necessarily unique)
 such that  $g_*^{m_i} = W(D_{x*}, D_{y*}^i)$.
 Let $n_i$ be the exponent sum of $D_{x*}$ in $W$.
 The numbers $m_i$ and $n_i$ can be computed on GAP,
 by using Reidemeister-Schreier
 style algorithms for subgroup presentations.
 
 The computation of $m_i$ is quite straightforward.
 The computation of $W_i$ (and hence $n_i$) requires
 a slightly more complicated, but standard, procedure.
 The idea is to have GAP compute a presentation for
 the subgroup of $SL_2(\mathbb{Z})$ generated by
 $D_{x*}$, $D_{y*}^i$ and $g_*^{n_i}$.
 After simplifying, GAP finds that the generator $g_*^{n_i}$ is redundant,
 and returns the word $W$. Details can be found in the source code at
 www.math.buffalo.edu/~jdmaster.

\section{Computations: slopes}

 Given an arbitrary $f \in \mathcal{M}_1^1$,
 we may associate slopes $\beta_3^1, \beta_3^2$ and $\beta_4^1, \beta_4^2$
 for the boundary of a cyclic cover $M_{f^m}$, as follows.
 
 To compute the slopes $\b_i^1$,
 We first compute an integer $m$ such that
 $f^m \in H_3 \cap H_4$.  We then compute words
 $W_3$ and $W_4$ such that $f_*^m = W_i(D_{x*}, D_{y*}^i)$,
 as described in the previous section,
 and let $n_i$ be the exponent sum
 of $D_{x*}$ in $W_i$.
  We then let $g_i = W_i(D_x, D_y^i) \in \mathcal{M}_1^1$.
 Since $g_i \in J_i$, then by the proof of Theorem \ref{main},
 we see that associated to the bundle $M_{g_i}$
 are slopes $\b_3^{1\prime} = (3, n_3)/gcd(3,n_3)$,
 $\b_4^{1\prime} = (2, n_4)/gcd(2,n_4)$,
 and $\b_4^{2 \prime} = (1,0)$.
  
 We use Lemma \ref{twist} to compute $t(f^m g_i^{-1})$, and then
 use Lemma \ref{frame} to compute that

\begin{eqnarray*}
 \beta_3^1 &=& (3, n_3 - 3(t(f^m g_3^{-1})))/gcd(3, n_3),\\
 \beta_4^1 &=& (2, n_4 - 2(t(f^m g_4^{-1})))/gcd(2, n_4),\\
 \b_4^2 &=& (1, -t(f^m g_4^{-1})).
\end{eqnarray*}

To compute $\b_3^2$, we first use Lemmas \ref{alex1}
 and \ref{alex2} to compute the Alexander polynomial
 for the relevant cover of the manifold $M_{g_3}$.
 Then, as in the proofs of Corollaries \ref{nonfiber1} and \ref{nonfiber2},
 we obtain the slope $\b_3^{2\prime}$.  Finally, using
 Lemma \ref{frame}, we compute the slope $\b_3^2$
 for the manifold $M_{f^m}$.

 If $\{\beta_3^1, \beta_3^2 \} \cap \{ \beta_4^1, \beta_4^2\} = \emptyset$,
 then we report ``success'',
 meaning that $M_f$ has a finite cover satisfying property VH.

\section{Example} Let $f = D_x^{-1} D_y$, so $M_f$ is the figure-eight
 knot exterior.  We shall show that $f^{12}$ has property VH.
 Using GAP, we compute:
  $m= 12$, so $f^{12} \in H_3 \cap H_4$;
 also $W_3(a,b) = (a^{-1} b a b a^{-1} b a b)^3$
 and $W_4(a,b) = (a^{-2} b^{-1} a^{-1} b^{-1} a^{-1})^4$,
 and so $n_3 = 0$ and $n_4 = -16$.
 We have $f_*^{12} =  W_3(D_{x*}, D_{y*}^3) = W_4(D_{x*}, D_{y*}^4)$.
 Letting $g_i = W_i(D_x, D_y^i)$,
 then  the slopes associated to $M_{g_i}$
  are $\b_3^{1 \prime} = (1,0)$, $\b_4^{1 \prime} = (1,-8)$,
 and $\b_4^{2 \prime} = (1,0)$.

 The Alexander polynomial, $\D_3(s)$, for the relevant cover
 of $M_{g_3}$ is given by:
\begin{eqnarray*} 
\D_3(s) &=& qs^4 + 3qs^3 +2qs^2 +2q s - q + q s^{-1} - 2q s^{-2}
 - 2q s^{-3} - 3q s^{-4} - q s^{-5}.
\end{eqnarray*}
 This polynomial is non-monic whenever $|q| > 1$. 
 Thus, associated to $M_{g_3}$ is the slope
 $\b_3^{2 \prime} = (1,0)$. 
 We use Lemma \ref{twist} to compute that $t(f^{12} g_3^{-1}) = -3$ and
 $t(f^{12} g_4^{-1}) = 4$.
 
 Then
\begin{eqnarray*}
 \beta_3^1 &=& (3, n_3 - 3(t(f^m g_3^{-1})))/gcd(3, n_3),\\
 &=& (1,3)\\
 \beta_3^2 &=& (1, -t(f^m g_3^{-1}))\\
 &=& (1, 3)\\
 \beta_4^1 &=& (2, n_4 - 2(t(f^m g_4^{-1})))/gcd(2, n_4),\\
 &=& (1, -12)\\
 \beta_4^2 &=& (1, -t(f^m g_4^{-1}))\\
  &=& (1, -4)
\end{eqnarray*}

 Since $\{\beta_3^1, \beta_3^2\} \cap  \{ \beta_4^1, \b_4^2 \} = \emptyset$,
 then $M_{f^{12}}$ has property VH.

\section{ An infinite family with Property VH}
\begin{proof} (of Theorem \ref{inf})

Let $f_n = (D_x^{-1} D_y)^{12} D_y^{12 n} $.  It may easily be checked
 that  the induced map on $H_1(F)$ has trace
 bigger than 3, and so the corresponding bundle
 $M_{f_n}$ is atoroidal.
 For $f_0$, a computation (see previous section) gives
 $\beta_3^1 =  (1, 3), \beta_3^2 = (1,3),
 \beta_4^1 = (1,-12)$ and $\b_4^2 = (1,-4)$.

 The element $f_n$ is equivalent in $SL_2(\mathbb{Z})$ to the
 element
 $$g_n = (D_x^{-1} D_y^3 D_x D_y^3 D_x^{-1} D_y^3 D_x D_y^3)^3 D_y^{12n}.$$
 One checks that the Alexander polynomial for the relevant cover
 of $M_{g_n}$
 is equivalent mod $n$ to the Alexander polynomial for the relevant
 cover of the bundle with monodromy
 $D_x^{-1} D_y^3 D_x D_y^3 D_x^{-1} D_y^3 D_x D_y^3$.
 One computes that that the leading coefficient of the latter polynomial
 is $q$, and thus the leading coefficient of the
 former polynomial is divisible by $q$.  Thus for the manifold $M_{g_n}$
 we have slopes $\b_3^{1 \prime} =  \b_3^{2 \prime} =(1,0)$.
 One computes $t(g_n f_n^{-1}) = -3$, and so
 for $f_n$, we have $\beta_3^1 = \b_3^2 = (1,3)$.

 Also, $f_n$ is equivalent in $SL_2(\mathbb{Z})$ to the element
 $h_n =  (D_x^{-2} D_y^{-4} D_x^{-1} D_y^{-4} D_x^{-1})^4 D_y^{12n}$.
 For the manifold $M_{h_n}$ we have $\b_4^{2 \prime} = (1,0)$,
 and we compute
 $\b_4^{1 \prime} = (1,-8)$.
 We compute $t(h_n f_n^{-1}) = 4$, and so
 for $f_n$, we have $\beta_4^1 = (1,-12)$, and $\b_4^2 = (1,-4)$.

 Since $\{ \beta_3^1, \beta_3^2 \} \cap \{ \beta_4^1, \b_4^2 \}
 = \emptyset$, then for all $n>0$, the manifold
 $M_{f_n}$ is finitely covered by a bundle with property VH,
 by Theorem \ref{main}.
 Furthermore, the manifolds $M_{f_n}$ are all
  obtained by doing surgery on the same hyperbolic
 knot $K \subset M_{f_0}$.
 Therefore, by \cite{LR} Section 3, there are infinitely many
 non-commensurable manifolds in the family $\{M_{f_n} \}$.
\end{proof}

\section{Computer results}

For every monodromy of complexity at most 5,
 the computer verified that a ``success'' criterion was met,
 and so the associated bundle is commensurable to one with property VH.
 The data for monodromies of complexity at most 5 is given below. 
 Since we are only interested in bundles up to commensurability,
 we have left out monodromies which are proper powers and monodromies
 with negative trace.
 Also, we have only included words up to cyclic permutations.
 The ``n'' in column two is the smallest positive integer
 such that $f^n \in H_3 \cap H_4$.
\\
\\
\begin{tabular} { l | l | l | l | l | l | l}
\underbar{monodromy f} & 
\underbar{n} &
\underbar{$t(f^n g_3^{-1})$} &
\underbar{$[\beta_3^1, \beta_3^2]$} &
\underbar{$t(f^n g_4^{-1})$} &
\underbar{$[\beta_4^1, \beta_4^2]$}&
\underbar{$M_{f^n}$ has VH?}
\\
&&&&&\\
$D_x^{-1} D_y$ & 12 & -3 & [(1,3),(1,3)]  & 4 & [(1, -12),(1,-4)]& Yes\\
$D_x^{-2} D_y$ & 12 & 2  & [(1,-6), (1,-2)] & 3 & [(1, -15),(1,-3)] & Yes\\
$D_x^{-1} D_y^2$ & 12 & -2 & [(1,6),(1,2)] & -3 & [(1,3),(1,3)] & Yes\\
$D_x^{-3} D_y$ & 6  & 2  & [(1,-6), (1,-2)] & 1 & [(1, -5),(1,-1)] & Yes\\
$D_x^{-2} D_y^2$ & 4 & -1 & [(1, 1),(1,1)] & -2 & [(1, 2),(1,2)] & Yes\\
$D_x^{-1} D_y^3$ & 6 & 0 &  [(1, -2),(1,0)] & -1 & [(1, 5),(1,1)] & Yes\\
$D_x^{-4} D_y$ & 4 & 1 &   [(1, -5 ), (1,-1)] & -2 & [(1, 2),(1,2)] & Yes\\
$D_x^{-3} D_y^2$& 12& 0 & [(1, -4),(1,0)] & -9 & [(1, 9),(1,9)] & Yes\\
$D_x^{-2} D_y^3$& 4 & 0 & [( 3, -8),(1,0)] & -1 & [(1, 1),(1,1)] & Yes\\
$D_x^{-1} D_y^4$& 4 & 1 & [(1, -3),(1,-3)] & 0 & [(1, -2),(1,0)] & Yes\\
$D_x^{-2} D_y D_x^{-1} D_y$ & 4 & 0 & [(3, -10),(1,-1)] & 0 & [(1, 0),(1,0)] &Yes\\
$D_x^{-1} D_y D_x^{-1} D_y^2$ & 12 & -6 & [(1, 10),(1,6)] & 0 & [(1, 0),(1,0)] & Yes
\end{tabular}
\\
\\
 From the table we see that
 every bundle of complexity at most five
 is commensurable to a bundle with Property VH.

 We considered all cyclically reduced, primitive, positive
 words on $D_x^{-1}$ and $D_y$
 of length at most 12, representing hyperbolic
 monodromies.  There are 745 of these.
 We verified that all but 36 of the associated bundles
 are finitely covered by a bundle with Property VH.
 Here we have not distinguished conjugate classes in
 $SL_2(\mathbb{Z})$, so
 these words do not all correspond to distinct bundles.
 The monodromy of smallest complexity which we cannot
 handle is $D_x^{-3} D_y^3$, for which $\b_3^1 = \b_4^1 = (1,-3)$.

 The routine runs quickly
 on words of rather large size. For example,  for
 the monodromy $f = D_x^{11} D_y^3 D_x D_y^6 D_x^{-4} (D_y D_x^{-1})^4 D_y$,
 a few seconds' computation gives
 $f^4 \in H_3 \cap H_4$, with
 associated slopes $\beta_3^1 = (1,-17),
 \beta_3^2 = (1,-15)$, and $\beta_4^1 = (1, -18), \b_4^2 = (1,-18)$.
 Thus $M_{f^4}$ has Property VH.

\vspace{5mm} \noindent Mathematics Department\\
 SUNY at Buffalo\\
Buffalo, NY 14260-2900\\
 jdmaster@buffalo.edu

\end{document}